\documentclass[11pt]{article}
\usepackage[pctex32]{graphics}
\usepackage{amssymb}
\usepackage{latexsym}
\usepackage{epsfig} 
\usepackage{multirow} 
\usepackage{pstricks}
\usepackage{float}
\usepackage{epstopdf}
\usepackage{mathrsfs}
\usepackage{soul}

\usepackage{amsmath}
\usepackage{amssymb}
\usepackage{relsize}
\usepackage{multirow}

\usepackage{graphicx}
\usepackage{multicol} 

\usepackage{times}
\usepackage{epstopdf}

\newtheorem{definition}{Definition}[section]
\newtheorem{theorem}[definition]{Theorem}

\definecolor{Red}{rgb}{1,0.,0.}

\usepackage{fancyhdr}
\usepackage{setspace}

\usepackage{hyperref}
\usepackage{mathtools}
\usepackage{algorithm}
\usepackage[dvipsnames]{xcolor}
\usepackage{algorithmic}
\usepackage{placeins}
\newcommand{\R}{{\mathbb R}}
\newcommand{\N}{{\mathbb N}}

\newcommand{\bbP}{{\mathbb P}}

\newcommand{\mT}{{\mathsf T}}
\newcommand{\mA}{{\mathsf A}}

\newcommand{\mR}{{\mathsf R}}
\newcommand{\mI}{{\mathsf I}}
\newcommand{\mC}{{\mathsf C}}
\newcommand{\mD}{{\mathsf D}}

\newcommand{\mSigma}{{\mathsf \Sigma}}
\newcommand{\mGamma}{{\mathsf \Gamma}}

\title{Discretization-free Bayesian inverse problems in distribution spaces}
\author{D Calvetti \and E Somersalo} 

\date{Case Western Reserve University\\
Department of Mathematics, Applied Mathematics and Statistics\\
Cleveland, OH, USA}

\begin{document}

\maketitle

\begin{abstract}
The Bayesian approach to inverse problems provides a practical way to solve ill-posed problems by augmenting the observation model with prior information. Due to its measure-theoretic underpinnings, the approach has raised theoretical interest, leading to a rather comprehensive description in infinite-dimensional function spaces. The goal of this article is to bridge the infinite-dimensional theory for linear inverse problems in distribution spaces and  associated computational inverse problems without resorting to a discrete approximation of the forward model. We show that the discretization of the unknown of interest is not necessary for the numerical treatment of the problem, the only approximations required being numerical quadratures that are independent of any discrete representation of the unknown. To demonstrate the viability of the approach, an analysis of X-ray tomography inverse problem is given in the proposed framework, and an analysis of the connection between the proposed approach and  a discretization-based one is also provided. 

\end{abstract}
\section{Introduction}

In most inverse problems the goal is to find a meaningful estimate of an unknown and unobservable quantity of interest based on noisy measurements of a related quantity that allows a direct observation. The bridge between the unknown of interest and the observed quantity is provided by a mathematical model,  often referred to as the forward or predictive model that may not have a well-defined or well-posed inverse. Thus the challenges of inverse problems arise from the ill-posed nature of inferring hidden causes from observed consequences. In the Bayesian interpretation, unknown quantities are modeled as random variables, the randomness representing epistemic uncertainty about their values as well as model uncertainties. A priori beliefs about the unknown of interest are encoded in the prior probability distribution, while the forward model with observation errors and model uncertainties constitute the likelihood model:  the posterior distribution obtained by Bayes' formula represent the solution to the inverse problem. 

Since the seminal work of Tarantola and Valette \cite{tarantola1982inverse} in the early 80s,  Bayesian methods have become one of the mainstream approaches to inverse problems \cite{tarantola2005inverse,kaipio2006statistical,calvetti2007introduction,calvetti2023bayesian}. The approach is particularly appealing because it naturally bridges classical inversion theory and  uncertainty quantification \cite{calvetti2018inverse}. While being fundamentally a practical approach to explicitly compute solutions in finite dimensional spaces, the measure theoretical underpinnings of the Bayesian approach have raised also significant interest in its formalism in the infinite dimensional function spaces where most physics-based models of inverse problems are originally formulated \cite{franklin1970well,mandelbaum1984linear,lehtinen1989linear,stuart2010inverse,dashti2015bayesian}.  A key question that has motivated much of the work in the infinite dimensional framework is whether and how the discretization, seen as an inevitable step for practical computations, affects the inverse problems theory. As stated by Stuart in \cite{stuart2010inverse}, the aim of the Hilbert space theory of inverse problems and its extensions to Banach spaces was to  {\em avoid discretization until the last possible moment.}  A similar motivation inspired \cite{lehtinen1989linear}, where the Bayesian theory for Gaussian measures was developed in distribution spaces. The cited works, and numerous follow-up articles contribute to our understanding of inverse problems in the infinite dimensional setting. In particular, one of the major motivations to develop the infinite dimensional theory of the Bayesian inversion is to understand the limit behavior of the analysis as the dimensionality of the problem grows, with the infinite dimensional model representing an ideal target of the limiting process. 
In particular, the analysis  has led to better converging sampling algorithms through the insight in the curse of dimensionality \cite{cotter2013mcmc}.
However, the gap between practical Bayesian computations and theoretically exact infinite-dimensional models prevails, even in the thoroughly understood realm of linear inverse problems with Gaussian prior and likelihood model. In fact, even if the infinite dimensional theory provides closed form solutions for the posterior mean and covariance operators in the function spaces, these formulas typically cannot be used in practical calculations without restating the problem in
 discretized and truncated form.

In this article, we revisit Gaussian linear inverse problems in infinite dimensional setting, modeling  unknown quantities of interest as generalized random variables, i.e., random variables taking on values in the space of distributions equipped with the weak topology and $\sigma$-algebra spanned by cylinder sets. Measurements are defined through dualities with test functions. We highlight that in this context {\em there is no need to discretize random variables}, by representing them in any basis spanning a function space. We postulate that while the observable is defined as a distribution, the data of inverse problems are a finite dimensional measurement of the observed random variable. The solution of the inverse problem is a probability distribution of a generalized random variable. To probe the random variable of interest, any finite set of test functions can be selected, and the probability distribution of the corresponding finite measurement can be computed. It will be shown that, under the Gaussian linear assumptions,  the computational  inverse problem can be completely solved without any discretization of the unknown of interest. Moreover, unlike in the basis-based Hilbert space theory, if the set of test functions probing the unknown of interest changes, the formulas for the mean and covariance of the posterior distribution make it possible to solve the modified problems without the need for restating the problem and carrying out the analysis from scratch. The passage of information from a finite measurement, i.e., the data, to any choice of finite  measurement of the unknown of interest occurs through a finite dimensional correlation matrix with entries that can be approximated numerically by quadratures, with no reference to any preselected basis of functions. It is in this sense that the theory developed here is discretization free, and devoid of any truncation errors due to basis selections.

The novelty of this article is twofold. On the conceptual level, we show that by framing the Bayesian inverse problems in distribution spaces, the finite dimensional theory can be used naturally to derive the formulas defining the posterior distributions without any reference to discretization of the unknowns (see Theorem 2.5). On the other hand, a practical example of X-ray tomography (see section~\ref{sec:tomography})  demonstrates that all calculations can be carried out without a reference to discretized variables, allowing to decide a posteriori which details of the unknown to extract from the posterior density, without the need to recompute anything, or to marginalize over details of no interest. As an added bonus, the approach provides new insight into the question of a consistent definition the tomography matrices that has been a topic of discussion in the literature.

The article is organized as follows. A review of the Bayesian framework for finite-dimensional inverse problems is given in section~\ref{sec:traditional}, where its extension to infinite-dimensional spaces is briefly discussed, identifying the challenges of this extension. In section ~\ref{sec:measurements} we set up the general framework of the proposed approach, interpreting observations as dual evaluations.  The derivation of the solution of the inverse problem as the posterior through dual evaluations is presented in subsection~\ref{sec:dual}. The formalism is elucidated by applying it to the inverse problem of X-ray tomography in section~\ref{sec:tomography}. In particular, in section~\ref{sec:comparison}, we highlight the differences between the proposed approach and the traditional one based on discretization.  For the sake of readability, the technical details are presented in the Appendix.

\section{Inverse problems in the Bayesian framework} \label{sec: BIP}

We start with a brief overview of the finite dimensional theory of linear inverse problems with Gaussian densities and a review of its extension to Hilbert spaces. We then develop the formalism for the current approach.

\subsection{Gaussian models in the traditional setting: a review}\label{sec:traditional}

We start by considering the linear finite dimensional inverse problem of estimating $x\in\R^n$ based on the observation 
\begin{equation}\label{lin model}
 b = \mA x + e,
\end{equation}
where $\mA\in\R^{m\times n}$ is a known matrix, $b\in\R^m$ is the observed datum, and $e\in\R^m$ represents observation noise. The Bayesian approach to inverse problems starts with the definition of a stochastic extension of the model (\ref{lin model}). 
Denoting by$(\Omega,{\frak S},{\mathbb P})$ a probability space, we write a stochastic model
\begin{equation}\label{lin model ext}
    B = \mA X + E,
\end{equation}
where the uppercase variables are defined as random variables, $X:\Omega\to \R^n$, and $B,E:\Omega\to\R^m$, the randomness reflecting the epistemic uncertainty about the values of the variables.
In a nutshell, the Bayesian solution to the inverse problem  (\ref{lin model})
through its stochastic extension (\ref{lin model ext}) is the posterior probability distribution of $X$ conditioned on $B=b$, where $b$ is the observed datum. For the sake of simplicity, assume that $X$ and $E$ are mutually independent, and their marginal probability distributions are absolutely continuous with respect to the Lebesgue measure in $\R^n$ and $\R^m$, respectively. The corresponding probability densities are denoted by $\pi_X$ and $\pi_E$, respectively.
Here, $\pi_X$ is referred to the prior density, encoding the level of uncertainty about $X$ prior to the measurement. Based on the mutual independency of $X$ and $E$, and the linear model (\ref{lin model ext}), the likelihood density of $B$, or density of $B$ conditional on $X=x$, is given by
\begin{equation}\label{lkh}
 \pi_{B\mid X}(b\mid x) = \pi_E(b - \mA x),
\end{equation}
i.e., the noise density is simply shifted around the presumably known value $\mA x$. According to Bayes' formula for probability densities, the posterior density of $X$, given the observation $B=b$, is proportional to the product of the prior and the likelihood densities,
\begin{equation}\label{Bayes}
\pi_{X\mid B}(x\mid b) \propto \pi_X(x)\pi_{B\mid X}(b\mid x),
\end{equation}
where ''$\propto$" stands for proportional up to a normalizing constant.

Of particular interest in this work are Gaussian random variables. Let $x_0\in\R^n$ and $\mGamma\in\R^{n\times n}$ be the mean and the covariance of the random variable $X$,
\[
 x_0 = {\mathbb E}(X) = \int_{\R^n} x \pi_X(x) dx, \quad
 \mGamma = {\mathbb E}\big(X-x_0)(X-x_0)^\mT\big) = 
 \int_{\R^n} (x-x_0)(x-x_0)^\mT \pi_X(x) dx,
\]
assuming that the integrals converge.  The random variable $X$ is Gaussian, or normally distributed, denoted by $X \sim {\mathcal N}(x_0,\mGamma)$, if its probability density is given by
\[
 \pi_X(x) \propto {\rm exp}\left(-\frac 12(x-x_0)^\mT \mGamma^{-1}(x-x_0)\right),
\]
and it is assumed that the symmetric matrix $\mGamma$ is positive definite.
If $X$ and $E$ are mutually independent and normally distributed,
\[
X \sim {\mathcal N}(x_0,\mGamma), \quad E\sim {\mathcal N}(0,\mSigma),
\]
where for simplicity we assumed vanishing mean of the noise, and the covariance matrices $\mGamma\in\R^{n\times n}$ and $\mSigma\in\R^{m\times m}$ are symmetric and positive definite (SPD), it is well-known that the posterior density is Gaussian,
\begin{equation}\label{Gaussians}
 X\mid (B=b) \sim {\mathcal N}(\overline x,\mC).
\end{equation}
The posterior mean and covariance matrix can be expressed in two equivalent ways as
\begin{eqnarray}
\overline x &=& \left(\mGamma^{-1} +\mA^\mT\mSigma^{-1}\mA  \right)^{-1}\big(\mGamma^{-1} x_0 + \mA^\mT\mSigma^{-1} b \big) =x_0+
\mGamma \mA^\mT \big( \mA \mGamma \mA^\mT + \mSigma \big)^{-1}(b - \mA x_0),\label{mean}\\
\mC &=& \left(\mGamma^{-1} +\mA^\mT \mSigma^{-1}\mA \right)^{-1} = \mGamma - \mGamma \mA^\mT\big( \mA \mGamma \mA^\mT + \mSigma \big)^{-1} \mA \mGamma,\label{cov}
\end{eqnarray}
the equivalence of the two alternative formulas being a direct consequence of the Sherman-Morrison-Woodbury formula. The above formulas are a direct consequence of the following theorem.

\begin{theorem}
Let $Z:\Omega\to\R^n$ and $Y:\Omega\to\R^m$ be two multivariate zero mean Gaussian random variables, and denote by $M$ the combined random variable,
\[
 M = \left[\begin{array}{c} Z \\ Y\end{array}\right] :\Omega \to \R^{n+m}, \quad M\sim{\mathcal N}(0, \mC), 
\]  
and partition joint covariance matrix as
 \[
 \mC = \left[\begin{array}{cc} \mC^{11} & \mC^{12} \\ \mC^{21} & \mC^{22}\end{array}\right] \in \R^{(n+m) \times (n+m)},
\]
where 
\[
 \mC^{11} = {\mathbb E}\big(Z Z^\mT), \quad \mC^{12} =  {\mathbb E}\big(Z Y^\mT) = \big(\mC^{21}\big)^\mT, \quad
\mC^{22} =  {\mathbb E}\big(Y Y^\mT).
\]
Then the posterior distribution of $Z$ conditioned on $Y=y$ is Gaussian with mean and covariance
\begin{eqnarray*}
\overline z &=& \mC^{12}\big(\mC^{22}\big)^{-1} y, \\
\mD &=& \mC^{11} - \mC^{12}\big(\mC^{22}\big)^{-1} \mC^{21} = \mC/\mC^{22},
\end{eqnarray*}
i.e., the conditional covariance is the Schur complement of $\mC^{22}$.
\end{theorem}

The formulas (\ref{mean}) and (\ref{cov}) follow by defining $Z = X-x_0$, $Y = B - \mA x_0$ and observing that
\[
 \mC^{11} = \mGamma, \quad \mC^{12} = \mGamma\mA^\mT, \quad  \mC^{22} = \mA\mGamma \mA^\mT + \mSigma.
\] 

The above linear Gaussian theory in finite dimensional spaces can be extended to the infinite-dimensional case. While the restriction to linear problems and Gaussian distributions is not needed and the theory covers more generally separable Banach spaces
\cite{stuart2010inverse,dashti2015bayesian}, we limit our discussion here to linear models in Hilbert spaces. For a general reference to Gaussian random variables in Hilbert spaces, see \cite{kukush2020gaussian,rozanov1971infinite}

Let $H_0$ and $H_1$ be two separable Hilbert spaces, equipped with the norms induced by the corresponding inner products,
\[
 \|x\|_0^2 =(x,x)_0, \quad \|y\|_1^2 = (y,y)_1,
\]
where $x\in H_0$, $y\in H_1$.
Let $X:\Omega\to H_1$ be a square integrable $H_0$-valued random variable,
\[
 {\mathbb E}\big(\|X\|^2_0\big) <\infty. 
\]
The square integrability guarantees that the mean $x_0\in H_0$  and covariance operator $\Gamma:H_0\to H_0$ of $X$,
\[
 {\mathbb E}(X) = x_0\in H_0,
\]
and
\[
{\mathbb E}\big((X-x_0,u)_0(X-x_0,v)_0\big) =(u,\Gamma v)_0, \quad u,v\in H_0,
\]
are well defined. 
We recall that $X$ is Gaussian with mean $x_0\in H_0$ and covariance $\Gamma:H_0\to H_0$, denoted by $X\sim{\mathcal N}(x_0,\Gamma)$, if for any finite collection of vectors $v_1,\ldots,v_n\in H_0$, the
multivariate random variable
\[
 X^n_c = \left[\begin{array}{c} (X-x_0,v_1)_0 \\ 
 \vdots  \\(X-x_0,v_n)_0\end{array}\right]:\Omega\to \R^n
\]
is Gaussian with zero mean and covariance matrix given by
\[
 {\mathbb E}\big((X-x_0,v_j)_0(X-x_0,v_k)_0\big) = \mGamma_{jk}^n = (v_j,\Gamma v_k)_0, \quad 1\leq j,k,\leq n.
\]

Assume that $A$ is a continuous linear map between two separable Hilbert spaces, $A: H_0\to H_1$. Further, let $X:\Omega\to H_0$ and $E:\Omega\to H_1$ be independent and Gaussian,
\[
X \sim {\mathcal N}(x_0,\Gamma), \quad E\sim {\mathcal N}(0,\Sigma).
\]
From the assumed square integrability of $X$ and $E$ it follows that the covariance operators $\Gamma:H_0\to H_0$ and $\Sigma:H_1\to H_1$ must be be positive semidefinite {\em nuclear}, or {\em trace class} operators: Given an orthonormal basis $\{v_k\}$ of $H_0$, we require that
\[
{\rm Trace}(\Gamma) = \sum_{j=1}^\infty ( v_j,\Gamma v_j)_0 <\infty,
\]
and similarly for $\Sigma$. 

Unlike in the finite dimensional theory, in the infinite-dimensional case, a reference to an underlying Lebesgue measure is meaningless, and therefore Bayes' theorem cannot be formulated in terms of densities as in (\ref{Bayes}). To obtain a proper extension of Bayes' theorem, we begin by denoting the prior Gaussian distribution by $\mu_X$ and the posterior density of $X$ conditional on $B=b$ by $\mu_{X\mid B}^b$. To state Bayes' theorem, we assume that
the posterior distribution is absolutely continuous with respect to the prior, and define the likelihood density of $B$, conditional on $X=x$ and denoted by $\pi_{B\mid X}(b\mid x)$, as the Radon-Nikodym derivative of the posterior density with respect to the prior, 
\[
 \pi_{B\mid X}(b\mid x) = \frac{d\mu^b_{X\mid B}}{d\mu_X}(x), \quad x\in H_0.
\]
This formula is the infinite dimensional equivalent of the Bayes' formula (\ref{Bayes}), as it implicitly defines the posterior density in terms of the prior and the likelihood.
Observe that in the case that $H_1$ is finite-dimensional, $H_1 = \R^m$ the likelihood density is explicitly given as
\begin{equation}\label{lkh2}
 \pi_{B\mid X}(b\mid x) \propto {\rm exp}\left( -\frac 12 \|b - A x\|^2_\mSigma\right),
\end{equation}
where the notation $\|z\|_\mSigma^2 = z^\mT \mSigma^{-1} z$ is used, and $\mSigma\in\R^{m\times m}$ is the matrix representation of the covariance operator in the canonical basis of $\R^m$. When the data space $H_1$ is infinite-dimensional, formula (\ref{lkh2}) is not immediately applicable, as $\Sigma$ is a trace class operator and therefore non-invertible. Therefore, in the infinite-dimensional setting, one can define the Cameron-Martin space $(C,\| \,\cdot\,\|_C)$ as a completion of the the space 
\[
 \{u\in H_1 \mid \| u\|_C = \| \Sigma^{-1/2} u\|<\infty\}.
\]
While in the finite-dimensional setting $E\subset H_1$ is the full space $H_1$, for infinite-dimensional spaces, by the Cameron-Martin theorem \cite{cameron1944transformations}, $\mu_E(C) = 0$, where $\mu_E$ is the Gaussian probability distribution of $E$. The functional 
\[
 \Phi(x) = \frac 12\|\Sigma^{-1/2}(b - A x)\|^2.
\]
 known as the Onsager-Machlup functional, has been studied in detail in the context of inverse problems in \cite{dashti2013map}.
One can show that formulas (\ref{mean}) and (\ref{cov}) have  corresponding equivalents in the infinite dimensional setting, however, the infinite-dimensional equivalent of the inverse of the matrix $\mC^{22}$ in the finite-dimensional model is no longer a bounded linear operator, but needs to be interpreted as a {\em measurable linear transformation}, see, e.g.,  \cite{mandelbaum1984linear,lehtinen1989linear}.  We refer to  \cite{stuart2010inverse,dashti2015bayesian} for further technical details of the Hilbert-space theory, and to \cite{lehtinen1989linear} for the extension to distribution spaces.

As pointed out before, a main motivation for developing the infinite-dimensional theory is to postpone the discretization to the last possible moment, i.e., to replace the standard ``discretize-then-analyze" approach by ``analyze-then-discretize" scheme. However, while it is reassuring to know that the finite dimensional theory has a consistent extension in the infinite-dimensional setting, the formulas for the infinite-dimensional random variables are not particularly useful in practice, as there is no immediate way of projecting the Gaussian posterior density to finite dimensional spaces. In practice, rather than discretizing the infinite-dimensional model, it is common to discretize and truncate the Hilbert space model and start the Bayesian analysis {\em de novo}, thus raising questions about the practical value of the analysis.

One motivation for the present contribution is to use the theory of distribution spaces while avoiding the complications in the infinite-dimensional theory, and to give a consistent finite dimensional description of the posterior density that bypasses the standard ``discretize-and-truncate"- process.

\subsection{Measurements as dual evaluations}\label{sec:measurements}

We start by outlining the general setting of this work. We denote by $U$ the topological vector space of test functions, e.g., the set of rapidly decreasing  $C^\infty$ test functions  in the Euclidean space, or compactly supported $C^\infty$ functions in $\R^n$, or periodic $C^\infty$-functions over the unit circle. 
We denote by $H$ the dual space of continuous linear functionals on $U$, consisting of the corresponding
distributions or generalized functions. 
The test function space is equipped with the Fr\'{e}chet topology, and it dual by the induced weak$^*$ topology.  
The duality between $U$ and $H$ is denoted by $\langle \varphi, x \rangle$, where $\varphi\in U$, $x\in H$.
Although topological considerations are not of central importance here, for the sake of making the article self-contained, some technical details are presented in the Appendix. For further details, see, e.g.,
\cite{rudin1991functional}.

We start with the definition of  measurement of a quantity in $H$.

\begin{definition}
Let $\{\varphi_1,\varphi_2,\ldots,\varphi_n\}$ be any finite collection of test functions  in $U$. 
A {\em measurement of a quantity $x\in H$ through the  measurement set $\Phi$} is a mapping,
\[
 H\to \R^n, \quad
 x\mapsto \langle\Phi, x\rangle  = \left[\begin{array}{c} \langle\varphi_1,x\rangle \\ 
  \langle\varphi_2,x\rangle\\
  \vdots\\
   \langle\varphi_n,x\rangle\end{array}\right], \quad
   \Phi  = \left[\begin{array}{c} \varphi_1 \\ 
  \varphi_2\\
  \vdots\\
   \varphi_n\end{array}\right] \in U^n.
 \]  
\end{definition}

The general theory of inverse problems developed in the framework of Bayesian computing \cite{kaipio2006statistical,calvetti2023bayesian} requires a stochastic extension of the deterministic model.  To this end, let $(\Omega,{\mathfrak S},\bbP)$ be a probability space, and define  {\em generalized random variables} \cite{gel2014generalized} as measurable mappings
\[
 X: \Omega \to H,
\] 
where $H$ is equipped with the $\sigma$-algebra induced by the cylinder sets of the form
\[
 S(u_1,\ldots,u_k) = \big\{x \in H\mid \big(\langle u_1,x \rangle,\ldots,\langle u_k,x \rangle\big) \in B, \; B\subset\R^k \mbox{ open, } u_1,\ldots,u_k\in U \big\},
\]
 that form the basis of the weak$^*$ topology of $H$.
 
 Of particular interest in this work are inverse problems arising from linear observation models. Let $U$ and $V$ be two test function spaces, and let $H = U^*$ and $K=V^*$ their topological duals. Consider a linear operator $A_0: V \to U$, and let $A = A_0^*:H\to K$ denote its adjoint extended to the duals: For a given $x\in H$, $A x\in K $  is defined through the identity
 \[
  \langle v ,A x\rangle = \langle A_0 v, x\rangle, \quad v \in V \mbox{ arbitrary.}
 \]
 For simplicity, we assume here that $A_0$ maps smooth functions to smooth functions. This assumption covers a large class of operators, such as pseudodifferential operators with smooth symbols that are pseudolocal, including, e.g., convolution operators with an integrable kernel. For a generalization to cases of non-smooth operators, we refer to the Appendix.
 In the classical deterministic setting, a linear inverse problem seeks to estimate a variable $x\in H$ based on the noisy observation of $b\in K$,
 \[
  b = A x + e,
 \]
 where $e\in K$ is a distribution representing the observation noise.  We extend the linear model to generalized random variables, following the notational convention to use uppercase letters to refer to random variables, and lowercase letters to their realizations. Thus the stochastic extension of the above equation can be written as
  \begin{equation}\label{forward}
  B = A X + E,
 \end{equation}
where $B$ and $E$ are $K$-valued generalized random variables, and $X$ is an $H$-valued random variable. We formalize the concept of a finite-dimensional measurement and the corresponding inverse problem formulation in the following definition.

\begin{definition}\label{def:IP} 
\begin{itemize}
\item[(a)] A finite-dimensional {\em direct measurement through a measurement set} $\big\{\varphi_1,\ldots,\varphi_n\big\}$  of a generalized random variable $X$ is defined as the random variable
\[
 M:\Omega\to\R^n, \quad M = \langle\Phi,X\rangle =  \left[\begin{array}{c} \langle \varphi_1,X \rangle \\ \vdots \\ \langle \varphi_n,X \rangle \end{array}\right].
\]  
\item[(b)] Given the forward model (\ref{forward}), the {\em indirect noisy measurement} of $X$ is any direct measurement
of the generalized random variable $B$,
\[
 Z :\Omega\to \R^m, \quad Z = \langle\Psi,B\rangle  =  \left[\begin{array}{c} \langle \psi_1,B \rangle \\ \vdots \\ \langle \psi_m, B\rangle \end{array}\right],
\]
where $\big\{\psi_1,\ldots,\psi_m\big\}\subset V$ is a given measurement set.  
\item[(c)]
The {\em generalized linear inverse problem} is to estimate any direct measurement $M = \langle\Phi,X\rangle$ of $X$ based on the given indirect measurement $Z = \langle\Psi, B\rangle$, where $X$ and $B$ are related to each other through the equation (\ref{forward}).  
\end{itemize}
\end{definition}

The above definition aims at keeping the approach practical in the sense that the data are finite, although defined in terms of an infinite dimensional model, and the solution $X$ is interrogated through arbitrary measurement sets, allowing the passage to an infinite dimensional limit, if necessary or desired.

In the following section, we analyze the problem in the Bayesian framework, and examine in more detail the case of Gaussian random variables, as the Gaussian setting is better suited for a detailed analysis of the differences between models based on discretization and the proposed discretization-free framework.

\subsection{Posterior density through dual evaluations}\label{sec:dual}

Consider the linear model (\ref{forward}), and a measurement of the quantity $B$ by a measurement set $\{\psi_k\}_{k=1}^m$. We have
\begin{equation}\label{observation j}
 \langle \psi_j,B\rangle =    \langle \psi_j,A X\rangle +  \langle \psi_j,E\rangle
   =    \langle A_0 \psi_j, X\rangle +  \langle \psi_j,E\rangle,
 \end{equation}
or
\begin{equation}\label{observation}
\langle \Psi,B\rangle = \langle A_0\Psi,X\rangle + \langle \Psi,E\rangle \
,
\end{equation}
where $A_0\Psi \in U^m$ corresponds to the measurement set $\{A_0\psi_j\}_{j=1}^m$.  In other words, an indirect noisy measurement with the measurement  $\Psi$ is equivalent to a direct noisy measurement with the measurement  $A_0\Psi$. Typically, for inverse problems, the measurement $A_0\Psi$ is not the preferred set unless we consider a simple denoising problem where $A_0$ is the identity. The practical formulation of the linear inverse problem in the Bayesian setting can be stated as follows.

\begin{definition}
Given the measurement sets $ \{\psi_k\}_{k=1}^m\subset V$ and $ \{\varphi_k\}_{k=1}^n\subset U$, find the posterior probability distribution of the $n$-variate random variable $\langle\Phi,X\rangle$ based on the observation $\langle\Psi,B\rangle$ given in (\ref{observation}).
\end{definition}

In a Hilbert space setting, the problem can be formulated in terms of projections: {\em Estimate the projection of $X$ on a finite dimensional subspace ${\rm span}\{\varphi_1,\ldots,\varphi_n\}$ based on the noisy observation of its projection on the subspace ${\rm span}\{A_0\psi_1,\ldots,A_0\psi_m\}$.} We point out that the problem is not requiring any discretization of the variables, but possibly only a numerical approximation for computing the projections.

We shall consider now the problem in the case of Gaussian distributions, further details being given in the Appendix. Let $X$ be an $H$-valued random variable. Recall that $X$ is a Gaussian generalized random variable if for every $n$ and every measurement set $ \{\varphi_j\}_{j=1}^n$, the $n$-variate random variables $\langle \Phi,X\rangle$ are Gaussian. The mean $m_X\in H$ and covariance operator $C_X:U\to H$ are defined by the identities
\[
 {\mathbb E}(\langle \varphi, X\rangle) = \langle \varphi, m_X\rangle,
\]  
\[
  {\mathbb E}\big((\langle \varphi, X\rangle -\langle \varphi, m_X\rangle)(\langle \phi, X\rangle -\langle \phi, m_X\rangle)
  \big) = \langle \varphi, C_X \phi \rangle,
\]  
where $\varphi,\phi\in U$ are arbitrary test functions. We use the standard notation $X \sim {\mathcal N}(m_X, C_X)$.  Similarly, we assume that the noise $E$ is a Gaussian $K$-valued generalized random variable, $E\sim{\mathcal N}(0, C_E)$, where $C_E:V\to K$. For simplicity, assume that $X$ and $E$ are mutually independent. Without loss of generality, we may assume that $X$ has zero mean, $m_X=0$, since adding the mean to the formulas afterwards is a straightforward matter.
The joint covariance matrix $\mC\in\R^{(n+m)\times (n+m)}$ of the pair $(\langle\Phi,X\rangle,\langle\Psi,B\rangle)$ is
\[
 \mC = \left[\begin{array}{cc} \mC^{11} & \mC^{12} \\ \mC^{21} & \mC^{22}\end{array}\right] \in \R^{(n+m) \times (n+m)},
\]
where
\begin{eqnarray}
\mC^{11}_{jk} &=& {\mathbb E}(\langle \varphi_j,X\rangle \langle \varphi_k,X\rangle  ) = \langle\varphi_j,C_X\varphi_k\rangle, \quad 1\leq j,k\leq n, \label{C11}\\
\mC^{12}_{jk} &=&  {\mathbb E}(\langle \varphi_j,X\rangle \langle \psi_k,A X + E \rangle  )  =  \langle\varphi_j,C_X A_0 \psi_k\rangle, \quad 1\leq j\leq n, \; 1\leq,k\leq m, \label{C12}\\
\mC^{21}_{kj} &=& \mC^{12}_{jk}, \\
\mC^{22}_{jk} &=&  {\mathbb E}(\langle \psi_j,A X + E\rangle \langle \psi_k,A X + E \rangle)  =
 \langle A_0 \psi_j,C_X A_0 \psi_k\rangle +  \langle \psi_j,C_E  \psi_k\rangle, \quad 1\leq j,k\leq m.\label{C22}
\end{eqnarray}
The posterior probability density is Gaussian,
\[
 \pi_{\langle\Phi,X\rangle\mid \langle\Psi,B\rangle}\big(z \mid\langle\Psi,b\rangle\big) = {\mathcal N}( z\mid \overline z,\mD),
\]
where   ${\mathcal N}( z\mid \overline z,\mD)$ refers to a Gaussian density with mean $\overline z$ and covariance matrix $\mD$, given by
\begin{eqnarray}
 \overline z &=& \mC^{12} \big(\mC^{22}\big)^{-1} \langle\Psi,b\rangle,  \label{mean2}\\
 \mD &=& \mC^{11} - \mC^{12}\big(\mC^{22}\big)^{-1} \mC^{21} = \mC/\mC^{22},\label{cov2}
\end{eqnarray} 
i.e., the posterior covariance is the Schur complement of $\mC^{22}$.

Consider now formula (\ref{mean2}).
Introducing the notation
\[
  \widetilde b = \big(\mC^{22}\big)^{-1}\langle\Psi,b\rangle,
\]
we write in component form as
\[
 \overline z_j  = \sum_{k=1}^m \langle\varphi_j,C_X A_0\psi_k\rangle \widetilde b_k = \left\langle\varphi_j,\sum_{k=1}^m C_X A_0\psi_k\widetilde b_k\right\rangle =
 \left\langle\varphi_j, C_X\big(A_0\Psi\big)^\mT \widetilde b \right\rangle,
\]
or in vector notation,
\begin{equation}\label{post mean}
 \overline z = \left\langle\Phi,
 C_X\big(A_0\Psi\big)^\mT \widetilde b \right\rangle,
\end{equation}
i.e., $\overline z$ is an evaluation through the measurement set $\{\varphi_1,\ldots,\varphi_n\}$ of the distribution
\[
 m_{X\mid(\Psi,B)} = C_X(A_0\psi)^\mT(\mC^{22})^{-1}\langle\psi,b\rangle \in H.
\]

We decompose the above mapping as
\[ 
\R^m\to H\to\R^n, \quad 
b\mapsto C_X(A_0\Psi)^\mT\,\widetilde b  \mapsto  \left\langle\Phi,
 C_X\big(A_0\Psi\big)^\mT\, \widetilde b \right\rangle.
\]
The first mapping represents the ``analysis" part, independent of any discretization of the unknown $X$, while the second step is the ``discretize" part. This organization genuinely postpones the discretization step to ``the last possible moment", which is the main motivation of the infinite-dimensional theory. We may write the computation of the posterior covariance in a similar fashion,
\[
 \mD = \left\langle\Phi,\left(\mI -C_X\big(A_0\Psi\big)^\mT (\mC^{22})^{-1} \big(A_0\Psi\big)\right)C_X\Phi\right\rangle,
\]
i.e., $\mD$ is the $n$-dimensional evaluation through the measurement set $\{\varphi_1,\ldots,\varphi_n\}$ of the covariance operator
\[
 C_{X\mid \langle\Psi,B\rangle}:U\to H, \quad
 C_{X\mid\langle\Psi,B\rangle} = C_X -C_X\big(A_0\Psi\big)^\mT (\mC^{22})^{-1} \big(A_0\Psi\big)C_X.
\]

We collect the results in the following theorem.

\begin{theorem}
Given a forward model (\ref{forward}), where $X$ and $E$ are Gaussian generalized random variables,
\[
 X\sim{\mathcal N(m_X,C_x)}, \quad E\sim{\mathcal N}(0,C_E),
\]
and a measurement set $\{\psi_1,\ldots,\psi_m\}\subset V$, we have
\[
 X\mid \langle\Psi,B\rangle \sim{\mathcal N}(m_{X\mid\langle\Psi,B\rangle},C_{X\mid\langle\Psi,B\rangle}),
\]
where the posterior mean and covariance operator are given by the formulas
\begin{eqnarray*}
m_{X\mid\langle\Psi,B\rangle}&=& m_X+ C_X(A_0\psi)^\mT(\mC^{22})^{-1}\langle\psi,b\rangle,\\
C_{X\mid\langle\Psi,B\rangle} &=&  C_X -C_X\big(A_0\Psi\big)^\mT (\mC^{22})^{-1} \big(A_0\Psi\big)C_X. 
\end{eqnarray*}
For any measurement set $\{\varphi_1,\ldots,\varphi_n\}\subset U$, we have
\[
\langle\Phi,X\rangle\mid\langle\Psi,B\rangle \sim {\mathcal N}(\overline z,\mD),
\]
where the finite-dimensional evaluations of the mean and covariance are obtained as
\begin{eqnarray*}
 \overline z &=& \langle\Phi,m_{X\mid\langle\Psi,B\rangle}\rangle, \\
 \mD &=& \langle\Phi,C_{X\mid\langle\Psi,B\rangle} \Phi\rangle.
\end{eqnarray*}
\end{theorem}
The gist of this theorem is that the mean and covariance operators are independent of any evaluation $\Phi$ of $X$, as well as of any discretization of the unknown of interest $X$. When the evaluation set $\Phi$ is chosen, the finite dimensional mean and covariance matrices are obtained through simple evaluations of distributions, without a need to recalculate the quantities if the evaluation set is changed.

In order to elucidate the challenges of actually computing the quantities, we work out the details for a standard inverse problem arising in X-ray tomography. 

\section{Example: X-ray tomography}\label{sec:tomography}

\begin{figure}[h!]
\centerline{
\includegraphics[width=0.9\textwidth]{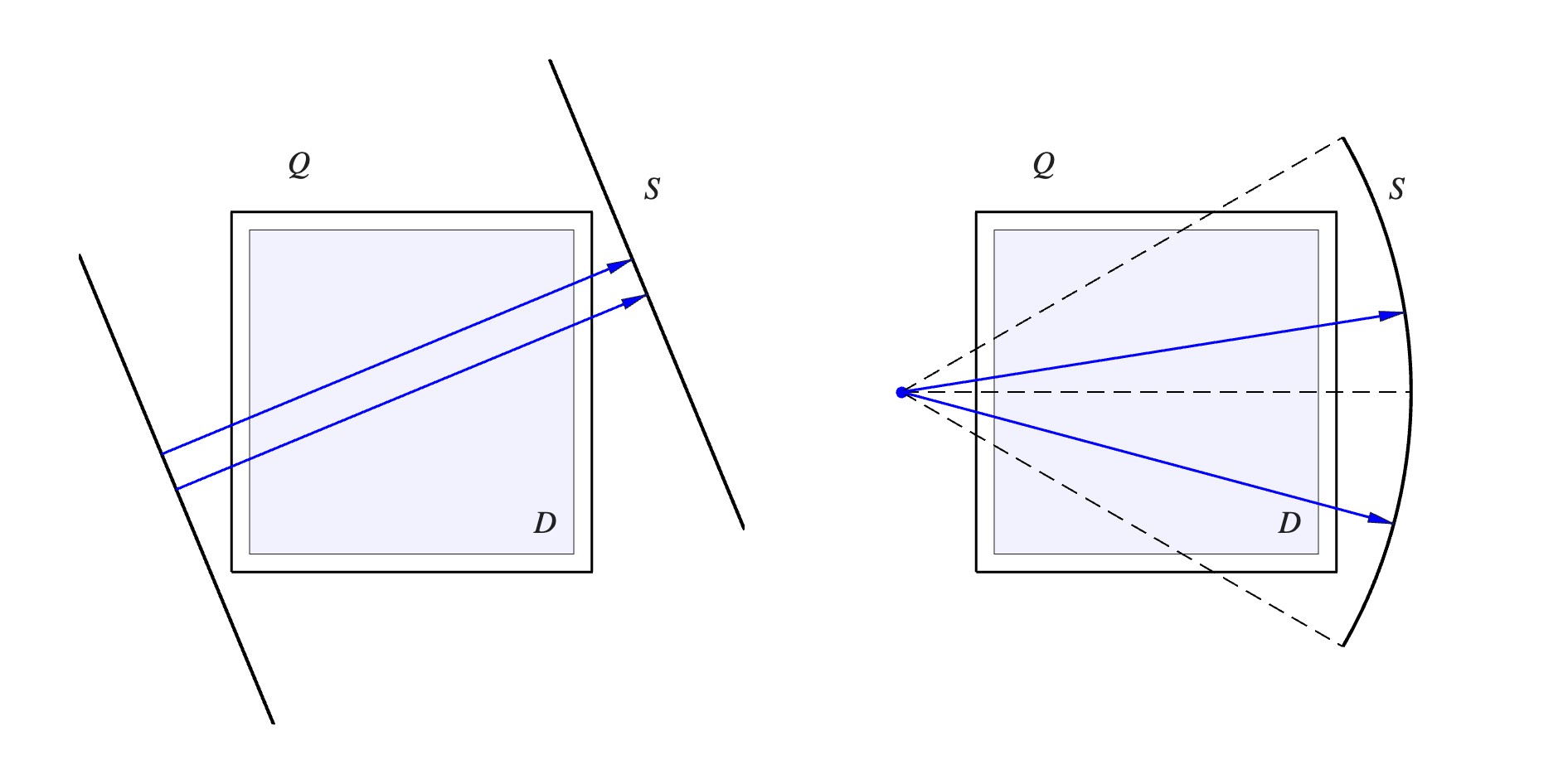}
}
\caption{ The geometric arrangement of a single projection in the parallel beam tomography (left) and the fan beam tomography measurement (right). The imaging window $Q$ containing the target is a square, and the unknown density is supported on the shaded square $D$. The source is either a line source (left) or a point source (right), and $S$ is the detector 
screen. The screen is assumed to  consist of non-overlapping intervals, each representing a single detectors collecting the photons falling on that detector and integrating the counts to a single pixel reading. To model this process, the measurement functions $\psi_k$, supported on the screen, are chosen to describe the detector sensitivities over those intervals. 
The full data consists of several projections obtained by rotating the device over a discrete set of projection angles around the target $Q$.}\label{fig: single fan}
\end{figure}

Consider the following measurement setting: Let $Q\subset\R^2$ denote the X-ray tomography image window containing an object of unknown density, and assume that $F$ represents the unknown density distribution in $Q$. We assume that $F$ is supported on a compact subset $\overline D\subset Q$, that is, $\langle \varphi, F\rangle = 0$ for every test function with ${\rm supp}(\varphi)\cap D = \emptyset$.  Typically, the domain $D$ is a rectangle.
The object is illuminated by X-rays emanating either from a line source (parallel beam tomography), or from a single point $p_0$ (fan beam tomography) outside the imaging window, and the attenuated radiation is measured along a screen $S$ at the opposite side of the object. In the case of the parallel beam tomography, $S$ is modeled as a line segment, while in the fan beam modality, we assume that the screen $S$ can be modeled by a curve segment on a circle of radius $R$ centered at the source $p_0$. The full tomographic data is collected by letting the measurement device rotate with respect to the target $Q$, yielding the projection data ideally from all directions.
A schematic description of these two geometries is shown in Figure~\ref{fig: single fan}.

For a given fixed illumination direction, we define the test function spaces as
\[
 U = C^\infty(Q), \quad V = C^\infty_0(S),
\]
with the corresponding dual spaces being
\[
 H = {\mathscr E}'(Q), \quad K = {\mathscr D}'(S),
\]
i.e., the distributions in $H$ are compactly supported on $Q$,  while the distributions in $K$ can be interpreted as restrictions on $S$ of one-dimensional distributions on the continuation of the line segment $S$.

\subsection{Forward model}

We derive the forward model for the distributions both for the parallel beam tomography and the cone beam tomography.

Starting with the parallel beam tomography, without loss of generality, let target region be defined as
$Q = [-1/2,1/2]\times [-1/2,1/2]$, and
assume that the projection direction is fixed. We parametrize the rays starting from the source line and finishing at the detector line as
\begin{equation}\label{ray}
 p(s,t) = s\omega^\perp(\theta) + t\omega(\theta),
\end{equation}
where
\[
 \omega(\theta) = (\cos\theta,\sin\theta), \quad \omega^\perp(\theta) = (-\sin\theta,\cos\theta),
\]
Here, $\theta$ is a fixed angle defining the detector orientation, and $t\in[-T/2,T/2]$ is the parameter of integration along the ray of length $T$ equal to the distance between the cource and the detector, while $s\in[-L/2,L/2]$ defines the source position along the line source of length $L$.

\begin{figure}[h!]
\centerline{
\includegraphics[width=10cm]{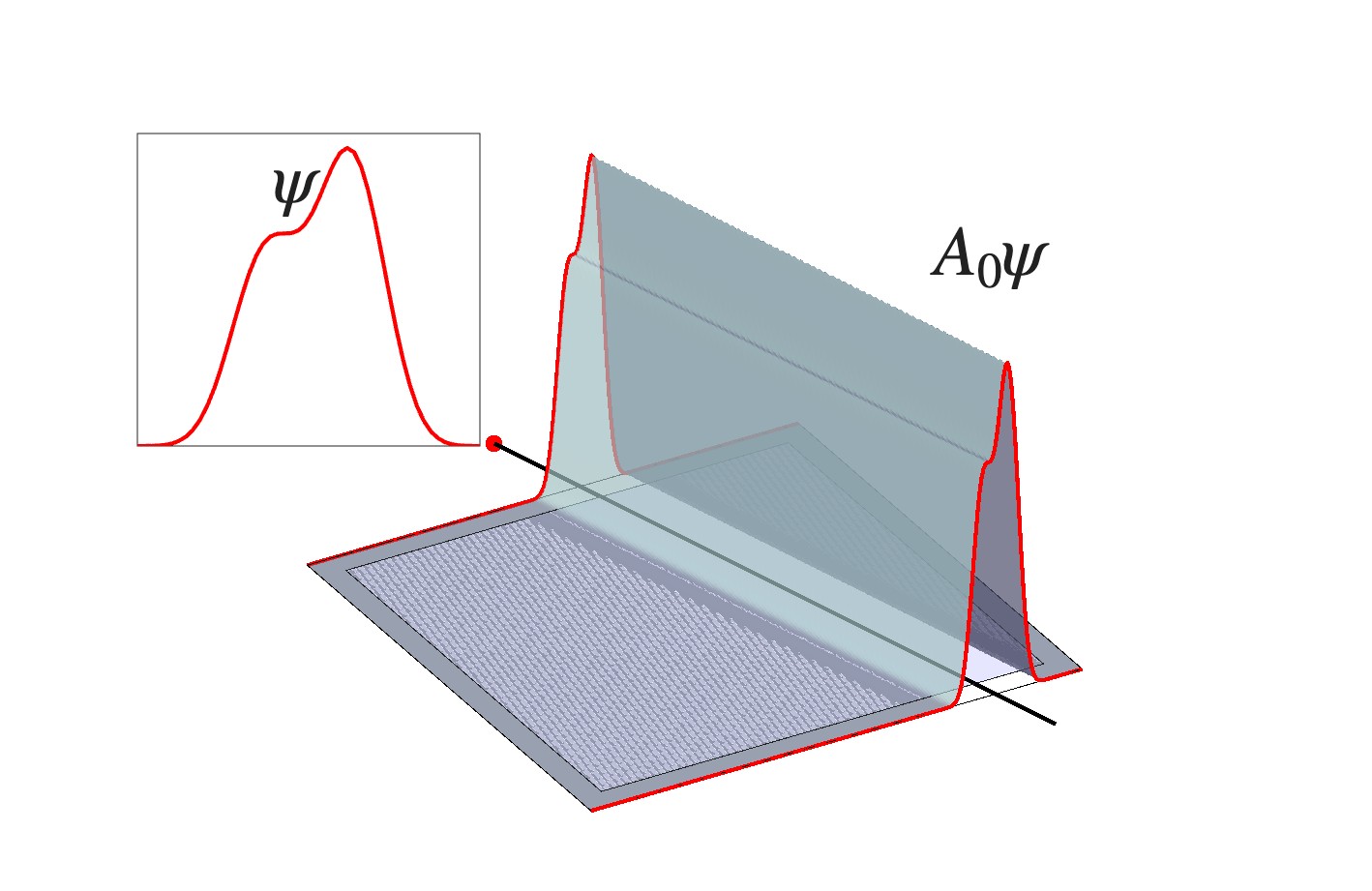}
\includegraphics[width=10cm]{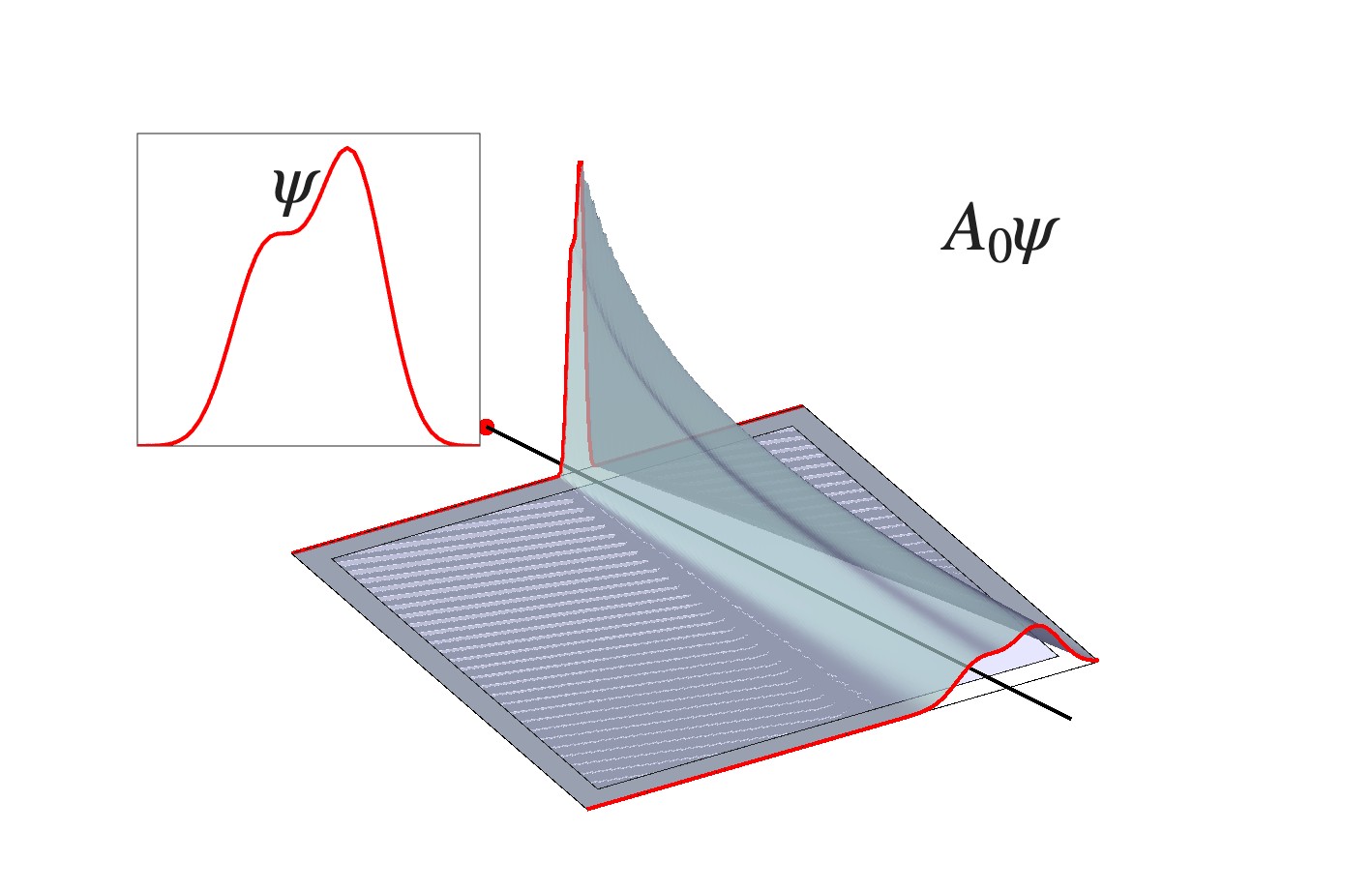}
}
\caption{The linear mapping $A_0:V\to U$ corresponding to a single parallel beam (left) and fan beam projection (right). The test function $\psi$ is a one--dimensional profile defined over the screen $S$, and $A_0\psi$ is a test function defined over the imaging area $Q$. Observe that in this visualization of the mapping $A_0$, the test function $\psi$ is an arbitrary compactly supported smooth function defined on the screen $S$,
not necessarily associated to any physical measurement device.}\label{fig:generalized tomo}
\end{figure}

Assuming first that $F=F(p)$ is a continuous density distribution of the material constituting the object, $p\in Q$, supported on $\overline D$, the line integrals along the ray (\ref{ray}) 
is given by
\[
 \int_{-T/2}^{T/2} F\big(
 s\omega^\perp(\theta) + t\omega(\theta)\big) dt.
\]
We assume that the screen $S = S_\theta$ consists of $m_\theta$ CCD detectors, modeled as disjoint intervals $s_k\subset S_\theta$, $1\leq k\leq m_\theta$, and that each detector integrates all photons reaching it into a single value. If the density function supported on $s_k$ is $\psi_k\in C^\infty_0(S_\theta)$, the $k$th data entry is given by
\[
 b_k(\theta) =\int_{-L/2}^{L/2} \psi_k(s) \int_{-T/2}^{T/2} F\big(
 s\omega^\perp(\theta) + t\omega(\theta)\big) dt ds.
\]
Observing that the integral extends over the rectangle
$Q' = [-T/2,T/2]\times[-L/2,L/2] $ and assuming that $D\subset Q'$, i.e., the rays fully cover the target, we may write the integral simply as
\[
 b_k(\theta) = \int_D \psi_k\big((\omega(\theta)^\perp)^\mT p\big) F(p) dp =\langle A_0\psi_k,F\rangle,
\]
where the forward operator $A_0$ is defined as
\begin{equation}\label{ridge function}
 A_0 : C^\infty_0(S_\theta) \to C^\infty(Q),\quad \psi_k \mapsto \widetilde 
 \psi_k,\quad
 \widetilde\psi_k(p) = \psi_k\big((\omega^\perp(\theta))^\mT p\big),
\end{equation}
i.e., the function $\widetilde\psi_k$ is a ridge function that is constant along lines parallel to the vector $\omega(\theta)$; see Figure~\ref{fig:generalized tomo} for an illustration. The same construction is repeated for all projection directions $\theta = \theta_1, \ldots,\theta_K$, giving rise to the full tomography operator $A_0$,
\[
  A_0 : C^\infty_0(S) \to C^\infty(Q), \quad S = S_{\theta_1}\dot\cup S_{\theta_2} \cdots \dot\cup S_{\theta_K},
\]
i.e., $S$ is the disjoint union of the screens $S_{\theta_k}$.

To derive the forward model for the fan beam geometry, assume again that the unknown $F$ is given as a classical continuous function supported on $\overline D$.
The idealized measurement with a single source-receiver position consists of integral of the function $F$ along rays from $p_0$ to $s\in S$, see Figure~\ref{fig: single fan}. Without loss of generality, we may assume that $p_0$ and $S$ are positioned symmetrically with respect to the horizontal line as in the figure, so that the data along $S$ can be modeled as
\[
 B(\theta)  = \int_0^R F(p_0 + t \omega(\theta)) dt, \quad \omega(\theta) = (\cos\theta,\sin\theta), \quad -\alpha<\theta<\alpha,
\] 
$R>0$ being the distance of the point source from the screen $S$.
Assuming again that the screen $S = S_{p_0}$ consists of $m_{p_0}$ pixels, e.g., CCD detectors, modeled as disjoint intervals $s_k\subset S$, $1\leq k\leq m_{p_0}$. Each detector integrates all photons arriving in $s_k$ into the single value. We may model the measurement by pixel values
\[
 b_k = \int_{-\alpha}^\alpha B(\theta)\psi_k(\theta) d\theta, \quad {\rm supp}(\psi_k)\subset s_k,
\]
where the measurement function $\psi_k$ is parametrized by the opening angle of the fanbeam and models the sensitivity distribution of the $k$th detector.  We write this integral as
\begin{eqnarray*}
 b_k &=& \int_{-\alpha}^\alpha \psi_k(\theta)\left(\int_0^R \int_{{\mathbb S}^1}F(p_0 + t\omega(\theta'))\delta(\theta-\theta')d\theta'\right) d\theta \\
 &=& \int_{{\mathbb S}^1}\int_0^R  \frac{\psi_k(\theta')}{t} F(p_0 +t\omega(\theta')) t d\theta' dt \\
 &=& \int_Q \widetilde \psi_k(p) F(p) dp \\
 &=& \langle A_0\psi_k, F\rangle,
\end{eqnarray*}
where we parametrized the domain $Q$ using polar coordinates with origin at $p_0$, noticing that ${\rm supp}(F)\subset Q$, and defined 
\begin{equation}\label{get A}
 \widetilde \psi_k(p) = (A_0 \psi_k)(p) = \frac 1{|p-p_0|} \psi_k(\theta), \quad \theta = {\arcsin}\frac{(p-p_0)_2}{|p-p_0|}.
\end{equation} 
Formula (\ref{get A})
defines a linear mapping
\[
A_0:C^\infty_0(S_{p_0})\to C^\infty(Q), \quad \psi_k\mapsto\widetilde\psi_k,
\]
visualized in the left panel of Figure~\ref{fig:generalized tomo}. To generate the full tomography data with several illumination directions parametrized by the point source position $p_0 = p_{0,1},\ldots,p_{0,K}$, the full tomography operator is defined as
\[
  A_0 : C^\infty_0(S) \to C^\infty(Q), \quad S = S_{p_{0,1}}\dot\cup S_{p_{0,2}} \cdots \dot\cup S_{p_{0,K}}.
\]

\subsection{Numerical approximation of the integrals}

We work out the details of the approximations of the matrix entries involved in the analysis in two special cases.

\begin{enumerate}
\item In the first case, we assume that the covariance operator $C_X$ is a scaled identity in the domain $D$, that is,
\begin{equation}\label{white noise prior}
 \phi\mapsto C_X \phi = \gamma^2 \chi_D\phi \in{\mathscr E}'(Q),
\end{equation}
where $\chi_D$ is the characteristic function of $D$ and $\gamma^2>0$, that is,
\[
 \langle \psi, C_X\phi\rangle = \gamma^2\int_D \psi(p)\phi(p) dp, \quad \psi,\phi\in C^\infty(Q).
\]
\item  In the second case, we assume that the generalized random variable $X$ is defined as a stationary Gaussian process (GP) with a continuous convolution-type kernel function $K$,
\begin{equation}\label{GP}
 C_X\phi(p) = \chi_D(p)\int_D K(p-q)\phi(q) dq,
\end{equation}
that is,
\[
 \langle\psi,C_X\phi\rangle =
  \int_D\int_D K(p-q)\psi(p)\phi(q) dq.
\]
The regularity properties of the realizations of a Gaussian process in terms of the properties of the covariance kernel $K$ are well studied in the literature, see, e.g., \cite{da2026sample} for a recent contribution.  The results cover regularity in both classical function spaces with H\"{o}lder continuity as well as in Sobolev spaces. A common choice is the Mat\'{e}rn kernel,
\[
 K(p-q) = \frac{2^{1-\nu}}{\Gamma(\nu)}\left(\sqrt{2\nu}\frac{|p-q|}\lambda\right)^\nu K_\nu\left(\sqrt{2\nu}\frac{|p-q|}\lambda\right),
\]
where $K_\nu$ is the modified Bessel function of the second kind, $\nu>0$ is the smoothness parameter, and $\lambda>0$ is the correlation length parameter. The reference to $\nu$ as a smoothness parameter is justified, as it can be shown that the sample paths are in the H\"{o}lder spaces  $C^{k,\delta}(D)$ for all $k$, $\delta$ such that $k+\delta<\nu$,
or in the $L^2$-based Sobolev scale, the sample paths are in $H^{s}(D)$ for $s<\nu+n/2$,
see \cite{da2026sample} for details.

While the sample path regularity is important in traditional spatial statistics involving estimates of pointwise evaluations and Reproducing Kernel Hilbert Space (RKHS) techniques, in the current setting where distributions are evaluated through dual pairings, the results are of less importance. 
\end{enumerate}

We consider first the parallel beam setting.
For simplicity, we assume that the noise is scaled white noise, that is,
for $\psi, \phi \in C^\infty_0(S)$, where $S=S_{\theta_\ell}$ for some $\ell =1,2,\ldots,K$,
\[
 \langle \psi,C_E \phi\rangle = \sigma^2 \int_{-L/2}^{L/2} \psi(s)\phi(s) ds.
\]
We start by considering the approximation of the matrix entries $\mC_{jk}^{22}$. With the above assumptions, the evaluations of the noise term is a straightforward integral,
\begin{equation}\label{noise contribution}
  \langle \psi_j,C_E \psi_k\rangle = \sigma^2 \int_{-L/2}^{L/2} \psi_j(s)\psi_k(s) ds,
\end{equation}
which can be approximated, e.g., by using a standard Gaussian quadrature rule.
Consider now the duality involving the forward map $A_0$. Assuming the prior model (\ref{white noise prior}) and using the notation (\ref{get A})
we have
\begin{equation}\label{int1}
 \langle A_0\psi_j,C_X A_0\psi_k\rangle =\gamma^2
 \int_D \widetilde\psi_j(p)\widetilde\psi_k(p) dp,
\end{equation}
and the numerical approximation of this integral can be done with a numerical quadrature rule over $D$. Similarly, using the model (\ref{GP}), we arrive at the integral
\begin{equation}\label{int2}
 \langle A_0\psi_j,C_X A_0\psi_k\rangle =\int_D
 \int_D \widetilde\psi_j(p)K(p-q)\widetilde\psi_k(q) dp dq,
\end{equation}
which, again, can be evaluated by quadrature rules.

\begin{figure}[ht!]
    \centerline{
\includegraphics[width=\textwidth]{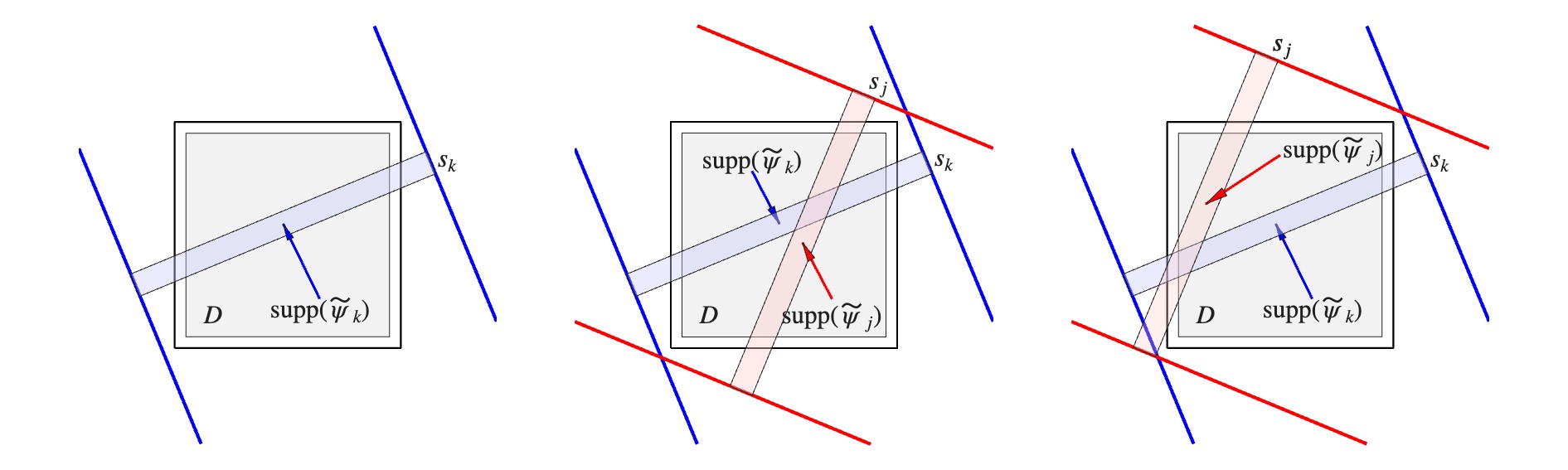} 
}
\caption{Computation of the matrix entries. Left: The integral (\ref{int1}) for $j=k$ and coinciding projection directions $\theta_\ell = \theta_{\ell'}$ equals the area of the parallelepiped $D\cap{\rm supp}(\widetilde\psi_k)$. Center: For $j$ and $k$ corresponding to different projection directions, the integral (\ref{int1}) equals the area of the parallepiped  
${\rm supp}(\widetilde\psi_j)\cap{\rm supp}(\widetilde\psi_k)$. Right: For $j$ and $k$ corresponding to different projection angles and supports of the ridge functions intersecting partly outside $D$, the integral (\ref{int1}) equals the area of the polygon
${\rm supp}(\widetilde\psi_j)\cap{\rm supp}(\widetilde\psi_k)\cap D$.
}\label{fig:matrix entries}
\end{figure}

The important observation here is that the formulas do not contain any reference to a ``tomography matrix". 
Furthermore, the integrals defining the matrix entries in some cases may be evaluated analytically. In particular, consider the case in which the CCD devices over the screen $S$ are modeled as $m$ identical disjoint subintervals covering the full screen, and furthermore, assume that the devices have a uniform sensitivity over the subinterval. We define
\[
 s_k = \left]\frac{k-1}{m}, \frac km \right[, \quad 1\leq k\leq m,
\]
and chose the test functions $\psi_k$ to be  $C^\infty$ approximations of the characteristic functions $\chi_{s_k}$ of the intervals $s_k$ to conform with the assumption of the uniform sensitivity.
With these assumptions, as the measurement functions approach the characteristic functions of the intervals $s_k$, the ridge functions $\widetilde \psi_k$ approach the characteristic functions of the strips of width $1/m$ across the imaging domains, and the integrals  (\ref{int1}) defining the matrix entries reduce to areas of polygonal domains.
More precisely, Let $\widetilde \psi_k = \widetilde\psi_{k,\theta_\ell}$ be a ridge function 
(\ref{ridge function}) with the projection angle $\theta_\ell$, and 
$\widetilde \psi_j = \widetilde\psi_{j,\theta_{\ell'}}$ be another one with projection angle $\theta_{\ell'}$, and assume that the measurement functions are  approximations of the characteristic functions of the detectors. Denoting by $\Gamma_k(\theta_\ell)$ the support of the ridge function (\ref{ridge function}), we have in the limit
\[
 \int_D \widetilde\psi_j(p)\widetilde\psi_k(p) dp = {\rm area}\big(\Gamma_j(\theta_{\ell'})\cap \Gamma_j(\theta_\ell)\cap D\big).
\]
In particular, for $\theta_\ell = \theta_{\ell'}$, the integrals vanish unless $k=j$. We refer to Figure~\ref{fig:matrix entries} for a geometric description of this formula.

A similar geometric reasoning for computing the matrix entries can be used also in the case of the fan beam tomography, however, even if the device functions $\psi_k$ are chosen as approximations of the characteristic functions over the detector, the integrals are slightly more complicated due to the presence of the reciprocal of the distance from the point source.

We turn now the attention to the computation of the cross correlation matrix entries. Given a set of measurement functions $\{\varphi_j\}$, the cross correlation matrix entries (\ref{C12}) are given by
\[
\mC^{12}_{jk} = \langle\varphi_j,C_X A_0\psi_k\rangle = 
\langle\varphi_j,C_X \widetilde\psi_k\rangle,
\]
and are represented by the integrals
\begin{equation}\label{cross corr1}
\mC^{12}_{jk} = \gamma^2 \int_D \varphi_j(p)\widetilde\psi_k(p) dp,
\end{equation}
if the prior model is the scaled white noise, or
\begin{equation}\label{cross corr2}
\mC^{12}_{jk} = \int_D \int_D \varphi_j(p)\widetilde\psi_k(q) K(p-q) dp dq,
\end{equation}
when a GP model is used as the prior.

As in the case of the matrix $\mC^{22}$, in particular cases, the cross correlation matrix can be evaluated exactly. Consider the choice of the measurement functions $\varphi_j$ to be $C^\infty$-approximations of characteristic functions of sets $P_j\subset D$, referred to as pixels. 
Consider the parallel beam measurement configuration with the device functions $\psi_k$ being approximations of the characteristic functions of the sets $s_k\subset S(\theta)$. In the limit when $\varphi_j$ converges to the characteristic function, the the integral (\ref{cross corr1}) converges to
\begin{equation}\label{R}
 C_{jk}^{12} = \gamma^2 {\rm area}\big(P_j \cap \Gamma_k(\theta)\big) = \gamma^2 \mR_{jk}(\theta), 
\end{equation}
where the matrix $\mR(\theta_\ell)$ can be identified with a standard tomography matrix as discussed in \cite{hansen2021computed,huber2025convergence,huber2025novel}. More generally, the computation of the cross correlation matrix entries reguire the use of numerical quadrature rules.

Before further discussion, we summarize here the main conclusions of this section. First, the discretization-free formulation of the inverse problem leads to relatively simple integral formulas for the matrix entries, and the numerical approximation of these entries can be done using numerical quadratures, while in some special cases, analytic formulas are available. Most notably, the calculation of the matrix $\mC^{22}$ requires no pixel representation of the unknown, nor a reference to a tomography matrix. As discussed in \cite{huber2025convergence,huber2025novel}, sometimes the tomography matrix and the matrix of the adjoint operator are incompatible in the sense that one is not a transpose of the other. In the present formulation, this becomes a moot point, as the tomography matrix appears only in the cross correlation operator.
The matrix $\mC^{22}$, as defined here, is automatically symmetric an positive definite, and 
the linear system involving this matrix can be solved, e.g., by using the conjugate gradient algorithm without the need to refer to the transpose of the tomography matrix. Second, the current formulation is in many ways economic, allowing to query only the part of the unknown of interest. More precisely, if one is interested in reconstructing the density image only in a subset $D'\subset D$, there is no need to model the unknown in the set $D\setminus D'$, but simply employ measurement functions $\varphi_j$  supported in $D'$. Effectively, the proposed approach generates a Bayesian estimate that in the traditional discretization-based formulation would need a marginalization process. The current approach could therefore be referred to as a marginalization-free approach. For a recent discussion of estimating the unknowns in subsets of the imaging domain, we refer to \cite{calvetti2025spotlight,calvetti2026spotlight}.

\subsection{Comparison with discretized model}\label{sec:comparison}

From the point of view of analyzing the potential gain of the proposed approach, it is useful to compare the formalism with the standard discretized model. To do so, assume that the imaging area $D$ is divided in pixels, denoted by $P_\ell\subset D$, $1\leq \ell\leq N$. For the sake of definiteness, we assume that the pixels are squares with non-intersecting interiors, and they cover the full imaging domain. To approximate the unknown $X$, we write
\[
 X \approx X^N = \sum_{\ell=1}^N X_\ell^N \chi_j,
\]
where $\chi_\ell$ is the characteristic function of $P_\ell$. Commonly, the univaraite random variables $X^N_\ell$ are identified with the pointwise value of $X$ at some interior point of $P_\ell$, if pointwise evaluation is meaningful. More generally, we may define $X^N_\ell$ through averaging over pixels
\[
 X_\ell^N =\frac{1}{|P_\ell|} \int_{P_\ell} X(p) dp,
\]
where $|P_\ell|$ is the area of $P_\ell$. The cross correlations of the pixel values are
\[
\mC^{(N)}_{\ell \ell'} =  {\mathbb E}\big(X_\ell^N X_{\ell'}^N\big) = \frac 1{|P_\ell||P_{\ell'}|}\int_{P_\ell}\int_{P_{\ell'}} {\mathbb E}\big(X(p) X(q)\big) dp dq, 
\]
which, in case of the white noise prior, reduces to
\begin{equation}\label{CN}
\mC^{(N)}_{\ell \ell'} = \frac{\gamma^2}{|P_\ell|} \delta_{\ell,\ell'},
\end{equation}
while in the GP prior case, we obtain
\[
\mC^{(N)}_{\ell \ell'} = \frac 1{|P_\ell||P_{\ell'}|}\int_{P_\ell}\int_{P_{\ell'}} K(p-q) dp dq. 
\]
Consider now the approximation of the matrix $\mC^{22}$. Using the discretized model, we obtain
\begin{eqnarray}
 {\mathbb E} \big(\langle A_0\psi_j,X\rangle\langle A_0\psi_k,X\rangle\big)
 &\approx&  {\mathbb E} \big(\langle A_0\psi_j,X^N\rangle\langle A_0\psi_k,X^N\rangle\big)\label{approx1}
 \\ &=&
 \sum_{\ell = 1}^N \sum_{\ell'=1}^N \mC_{\ell,\ell'}^{(N)} \int_{P_\ell} \widetilde\psi_j(p) dp \int_{P_{\ell'}} \widetilde\psi_k(q) dq.\label{approx2}
\end{eqnarray}
We generalize the formula (\ref{R}), and define the tomography matrix as
\begin{equation}\label{tomomat}
\mR_{\ell j} = \int_{P_\ell}\widetilde\psi_j(p) dp.
\end{equation}
Observe that in the spacial case of the parallel beam tomography when $\widetilde\psi_j$ is the ridge function corresponding to a characteristic function of a sensor at a given projection direction $\theta$, the definition of $\mR$ in (\ref{tomomat}) coincides with (\ref{R}).
With this notation, we arrive at the approximation
\[
 {\mathbb E} \big(\langle A_0\psi_j,X\rangle\langle A_0\psi_k,X\rangle\big) \approx \big(\mR^\mT \mC^{(N)} \mR\big)_{jk}.
\]
Assuming scaled white noise model for the measurement noise $E$, we therefore arrive, maybe not surprisingly, at the standard approximation,
\[
 \mC^{22} \approx \mR^\mT\mC^{(N)}\mR + \sigma^2 \mI.
\]
It is of interest to see how the discrete approximation is related to the approximation of integrals. For simplicity, consider the case of scaled white noise prior.
Using the model (\ref{CN}), the formula (\ref{approx2}) assumes the form
\[
 \sum_{\ell = 1}^N \sum_{\ell'=1}^N \mC_{\ell,\ell'}^{(N)} \int_{P_\ell} \widetilde\psi_j(p) dp \int_{P_{\ell'}} \widetilde\psi_k(q) dq
 ={\gamma^2} \sum_{\ell = 1}^N
 \frac 1{|P_\ell|}\int_{P_\ell}\widetilde\psi_j(p) dp 
\int_{P_\ell}\widetilde\psi_k(p) dp.
\]
Comparing with the discretization free form (\ref{int1}), we see that the pixel-based approximation is equivalent to 
\[
\int_D \widetilde\psi_j(p)\widetilde\psi_k(p) dp \approx
\sum_{\ell = 1}^N
\frac 1{|P_\ell|}\int_{P_\ell}\widetilde\psi_j(p) dp 
\int_{P_\ell}\widetilde\psi_k(q) dq.
\]
The error of this approximation can be written as
\begin{eqnarray*}
 \varepsilon^N_{jk} &=& 
\sum_{\ell=1}^N\int_{P_\ell} \widetilde\psi_j(p)\widetilde\psi_k(p) dp -
\sum_{\ell = 1}^N
\frac 1{|P_\ell|}\int_{P_\ell}\widetilde\psi_j(p) dp 
\int_{P_\ell}\widetilde\psi_k(q) dq \\
&=&\sum_{\ell = 1}^N \int_{P_\ell}\widetilde\psi_j(p)\left\{\widetilde\psi_k(p) - \frac 1{|P_\ell|}\int_{P_\ell}\widetilde\psi_k(q)dq\right\} dp,
\end{eqnarray*}
yielding an error estimate
\[
|\varepsilon^N_{jk}| \leq \max_{1\leq\ell\leq N}\sup_{p\in P_\ell}\left|
\widetilde\psi_k(p) - \frac 1{|P_\ell|}\int_{P_\ell}\widetilde\psi_k(q)dq\right|\int_D|\widetilde\psi_j(p)|dp,
\]
which converges to zero as $N$ increases by the uniform convergence of the mean value of the $C^\infty$-functions $\widetilde\psi_k$.

\section{Discussion}

The present article aims at bridging the computational Bayesian inverse problems methodology and the theoretical approaches in infinite dimensional Hilbert and Banach spaces in the framework of distribution spaces. The distributional environment provides an infinite-dimensional context for the theory, however, the evaluations of distributions through applications on test functions are finite-dimensional, and the finite dimensional evaluations are the basis for both the topology and the $\sigma$-algebras of these spaces. The evaluations through applications by test functions are also conceptually in line with the idea of measurements in engineering, e.g., by integrating incoming radiation over a light sensor, or integrating magnetic flux density over a magnetometer coil, or current density over an electrode, or averaging a signal over a time window, to mention a few examples, highlighting the convergence of theoretical and practical viewpoints.

In this article, it has been shown that to develop a consistent Bayesian framework for inverse problems in distribution spaces, one can define the forward model between infinite dimensional distribution spaces without certain technical challenges arising from the fact that covariance operators in Bananch spaces are compact operators. Moreover, it was demonstrated via the X-ray tomography example that the discretization of the unknown is in general not necessary for the numerical evaluation of the covariance matrices needed in the practical computations. Furthermore, in the Gaussian linear theory discussed in this article, once the inverse problem  has been solved in the distribution space, one can decide a posteriori which test functions are used to probe the posterior density without the need to recalculate the posterior mean and covariance, which are independent of the a posteriori discretization as Theorem 2.5 shows.
The framework therefore provides a natural and flexible framework, e.g., for querying detailed spatial information of the unknown without a need to calculate marginal distributions.

While the approach leads to useful and computationally applicable results, the aims of this paper are more conceptual than computational. For this reason, detailded computed examples are not included in this article, but are left for future investigations.
Furthermore, to underline the conceptual structure, here the forward model was assumed to arise from  an operator that maps smooth functions to smooth functions. As pointed out in the Appendix, this simplification is not necessary, and the functional analytic details depend on the particular problem at hand.  The extension of the proposed discretization-free approach to nonlinear problems, and  a discussion of how the design algorithms that interface well with matrix-free iterative linear solvers will be topics of future investigations.

\section*{Acknowledgements}

The work of DC was partly supported by the National Science Foundation grants DMS 1951446 and DMS-2513481, and that of  ES by the National Science Foundation grants DMS-2204618 and DMS-2513481. Support from the John Simons Guggenheim Foundation to ES and from Simons Foundation to DC is acknowledged with gratitude.

\section*{Appendix}

For the convenience of the reader, in this appendix we summarize some central results concerning the functional analytic framework on which the analysis is based. For further details, we refer to \cite{lehtinen1989linear} and references therein.

Let $U$ be a separable complete vector space equipped with a system of inner products $(\,\cdot\, ,\,\cdot \,)_n$, $n=1,2,3.\ldots$ such that the corresponding induced norms $\|\,\cdot\, \|_n$  are non-decreasing, i.e., 
$\|u\|_n\leq \|u\|_{n+1}$ for all $u\in U$ and $n=1,2,\ldots$. We denote by $U_n$ the completion of $U$ with respect to the norm $\|\,\cdot\,\|_n$, defining a separable Hilbert space. It follows from the monotonicity of the norms that $U_1\supseteq U_2 \supseteq \ldots$, and 
\begin{equation}\label{U}
 U = \bigcap_{n=1}^\infty U_n.
\end{equation}
The space $U$ defined as an intersection of nested Hilbert spaces (\ref{U}) is referred to as a {\em countably-Hilbert space}.  Its topology is spanned by neighborhoods of zero of the form $\{u\in U\mid \|u\|_n <\varepsilon\,\; 1\leq n\leq N\}$,  $N\geq 1$, $\varepsilon>0$.

Let $H$ denote the space of continuous linear functionals of $U$, equipped with the weak$^*$ topology, and let $H_n$ denote the dual space of $U_n$. 
We have $H_1\subseteq H_2\subseteq \ldots$, and we may identify $H$ as
\begin{equation}\label{H}
 H = \bigcup_{n=1}^\infty H_n,
 \end{equation}
 We denote the duality $U\times H\to \R$ (or ${\mathbb C}$) by $(u,\xi) \mapsto \langle u,\xi\rangle$.

A bounded linear operator $B:U_n\to U_n$ is said to be {\em nuclear}, or {\em trace-class}, if 
\[
 {\rm Trace}(B) = \sum_{j=1}^\infty(v_j,B v_j)_n,
\]
where $\{v_j\}_{j=1}^\infty$ is any orthonormal basis  of $U_n$. 
The countably-Hilbert space is called {\em nuclear}, if for every $k$,  the bilinear functional $B_k$ representing the $k$th inner product,  $B_k(u,v) = (u,v)_k$, $u,v\in U$, is represented in some $U_n$, $n>k$, as
\[
 B_k(u,v) = (u,B_k^n v)_n,\quad u,v\in U,
\]
where the operator $B_k^n:U_n\to U_n$ is nuclear.  

To make these concepts more concrete, consider periodic $C^\infty$-functions on the unit circle,
\[
 U = C^\infty({\mathbb S}^1),
\]
and the inner products 
\[
 (u,v)_n = \sum_{j=-\infty}^\infty \langle j \rangle^{2n} \widehat u_j \overline{\widehat v_j} ,
\]
where $\langle j\rangle = (1+j^2)^{1/2}$, and $\{\widehat u_j\}_{j=-\infty}^\infty$ is the Fourier transform of the function $u$. 
By Parseval's theorem, the completion of $U$ with respect to the norm $\|u\|_n = \big\|(\langle j\rangle^n \widehat u_j)_{j=-\infty}^\infty\big\|_{\ell^2}$ defines the Sobolev space with smoothness index $n$,  $W^n({\mathbb S}^1) = U_n$.
The norms are non-increasing, and for every  $k$ and $n>k$, we have
\[
  (u,v)_k = \sum_{j=-\infty}^\infty \langle j \rangle^{2k} \widehat u_j \overline{\widehat v_j}  =  \sum_{j=-\infty}^\infty \langle j \rangle^{2n} \big(\langle j\rangle^{2(k-n)} \widehat u_j \big)\overline{\widehat v_j} 
  =(B_k^n u,v)_n,
\]  
where $B_k^n$ is the Fourier multiplier 
\[
 B_k^n: u \mapsto \sum_{j=-\infty}^\infty e^{ij\theta} \langle j\rangle^{2(k-n)} \widehat u_j,
\] 
which is trace class operator for $n\geq k+1$, showing that $U$ is a nuclear countably-Hilbert space.  
We define the dual pairing through extension of the $L^2$-inner product.
In this example, the space $H$ consists of the distributions on the unit circle, $H = {\mathscr D}'({\mathbb S}^1)$, and the spaces $H_n$ can be identified with Sobolev spaces with negative smoothness index, $H_n = W^{-n}({\mathbb S}^1)$.

Another example relevant for the present discussion is the Schwartz class of rapidly decreasing $C^\infty$-functions,
\[
 U = {\mathscr S}(\R^n) =\big\{ u\in C^{\infty}(\R^n)\mid |x^\alpha D^\beta u(x)|\leq C_{\alpha,\beta} \mbox{ for all $\alpha,\beta \in \N_0^n$}\big\},
\] 
the space $H = {\mathscr S}'(\R^n)$ being the class of tempered distributions,
and 
\[
 U = {\mathscr D}(\Omega) = \big\{ u\in C^{\infty}(\Omega)\mid {\rm supp}(u) \mbox{ compact}\subset \Omega \big\}
 \]
 for some bounded smooth $\Omega\subset\R^n$, and $H = {\mathscr D}'(\Omega)$ are distributions over the set $\Omega$. The construction of the countably-Hilbert space structure for this example is along the lines of the discussion for the previous example, requiring more technical operator-theoretic details that   can be found, e.g.,  in \cite{reed2012methods} for the former, and \cite{rozanov1971infinite} for the latter.  
 
Given a probability space $(\Omega,{\frak S},{\mathbb P})$, consider a measurable mapping $\xi: U \to L^2(\Omega,{\frak S},{\mathbb P})$. If $\xi$ is continuous and linear, we call $\xi$ a {\em linear random functional} (LRF). We define the mean $\mu$ and correlation functional $C$ of $\xi$ by the formulas
\[
 \mu(u) = {\mathbb E}\big(\xi(u)\big), \quad C(u,v) = {\mathbb E}\big((\xi(u)-\mu(u))(\xi(v)-\mu(v)\big).
\] 
The continuity of the mapping $\xi$ implies that the mean and correlation functional admit representations of the form 
\[
 \mu(u) = (u,a_n)_n, \quad C(u,v) = (B_n u,v)_n
\]
for some $n\geq 1$, where $a_n\in U_n$ and $B_n\in L(U_n)$, i.e., $B_n$ is a linear bounded operator in $U_n$.  Consequently, the LRF can be extended to the Hilbert space $U_n$, thus defining an LRF in the Hilbert space
$U_n$.

A  LRF on $U$ is Gaussian if, for every $u_1,\ldots,u_k\in U$, the multivariate random variables $(\xi(u_1),\ldots,\xi(u_k))$ are Gaussian.  In the following, we restrict the discussion to Gaussian random variables.

To define random variables in the dual space $H$, we equip it with a $\sigma$-algebra spanned by the cylinder sets of the form 
\[
 {\mathcal C} = \big\{\xi\in H\mid  (\langle u_1,\xi\rangle,\ldots,\langle u_k,\xi\rangle) \in B,\; u_1,\ldots,u_k\in U, \; B\subset \R^k \mbox{ open} \big\}.
 \]
It turns out that this $\sigma$-algebra coincides with the Borel $\sigma$-algebra relative to the weak$^*$ topology, see \cite{rozanov1971infinite}.
A natural question is whether a LRF defines an $H$-valued {\em generalized random variable}
\begin{equation}\label{Xi}
\Xi:\Omega \to H \mbox{ such that $\xi(u) = \langle u, \Xi\rangle$ for all $u\in U$.}
\end{equation}
In the affirmative case, we refer to the random variable $\Xi$ as a {\em Gaussian linear functional}.
As demonstrated in \cite{rozanov1971infinite}, for a Gaussian LRF, the necessary and sufficient condition for the existence of the representation as a Gaussian linear functional is that the covariance operator $B_n$ is a trace-class operator with respect to some orthonormal basis of $U_n$. Since we have assumed that the countably-Hilbert space $U$ is nuclear, this condition is automatically satisfied, as $B_n$ can always be interpreted as a trace-class operator in $U_{n+k}$ for some $k>0$.  We therefore conclude that {\em in a nuclear countably-Hilbert space $U$, every Gaussian LRF can be represented as a Gaussian linear functional (\ref{Xi})}

Finally, we discuss briefly some extension of the formalism to cases where the forward model $A$ cannot be defined via an operator $A_0:V\to U$, implicitly assuming that $A_0$ maps smooth functions to smooth functions, a condition that is violated by, e.g., integral operators of non-convolutional form with a non-smooth kernel. Consider a linear mapping $A_0$ that maps the test functions in $V$ to a distribution, $A_0:V\to H$. Without restricting significantly the generality, assume that for some $n\geq 1$ we have ${\rm Im}(A_0) = \{A_0 v\mid v\in V\} \subset H_n$, and assume that $A_0:V\to H_n$ is continuous. Let $J_n: H_n \to U_n$ denote the canonical representation of the dual of the Hilbert space $H_n$ in  $U_n$, $\langle \xi,u\rangle = (J_n\xi,u)_n$,
and define $A_n = J_n A_0: V\to U_n$. We can then define the forward map as an adjoint $A =A_n^*:H_n\to K$, and the inverse problem is defined for a random variable $X$ assuming that $X\in H_n$ with probability one. For more detailed continuity considerations of the adjoint $A = A^*_n$, we refer to \cite{lehtinen1989linear}.

In \cite{lehtinen1989linear}, where the general theory of linear inverse problems for Gaussian generalized random variables was developed,  it was shown that the posterior mean and covariance operator are well-defined in a sense compatible with the earlier analysis in \cite{mandelbaum1984linear}. Since the present article does not lean on those interpretations, we will not go further into that theory here.

 \bibliographystyle{siam} 
\bibliography{ThOfMeas_biblio}

\end{document}